\documentclass[a4paper,11pt]{ltdartcl}
\usepackage[
    backend=biber,
    style=ieee,
    url=false,
  ]{biblatex}
\usepackage[T1]{fontenc}

\usepackage{amsmath}
\usepackage{mathtools}
\usepackage{amsthm}
\usepackage{amsfonts}  
\usepackage{amssymb}
\usepackage{mathrsfs}
\usepackage{xcolor}
\usepackage{upgreek}
\usepackage{hyperref}
\newcommand\R{\mathbb{R}}
\usepackage{graphicx}
\newtheorem{theorem}{Theorem}[section]

\title{Second order optimality conditions in a new Lagrangian formulation for optimal control problems}

\author
{Michael Konopik\footnote{Friedrich-Alexander-Universität Erlangen-Nürnberg (FAU), Institute of Applied Dynamics (LTD), Immerwahrstrasse 1, 91058 Erlangen, Germany. Email: \href{mailto:michael.konopik@fau.de}{michael.konopik@fau.de}}\ \thanks{The work of this author has been supported by Deutsche Forschungsgemeinschaft (DFG), Grant No. LE 1841/12-1, AOBJ: 692092.}
\qquad
Sigrid Leyendecker\footnote{Friedrich-Alexander-Universität Erlangen-Nürnberg (FAU), Institute of Applied Dynamics (LTD), Immerwahrstrasse 1, 91058 Erlangen, Germany. Email: \href{mailto:sigrid.leyendecker@fau.de}{sigrid.leyendecker@fau.de}}\ %\thanks{The work of this author has been supported by Deutsche Forschungsgemeinschaft (DFG), Grant No. LE 1841/12-1, AOBJ: 692092.}
\qquad 
 Sofya Maslovskaya\footnote{\textit{First author}, \textit{corresponding author},  Universität Paderborn (UPB), Numerical Mathematics and Control (NMC), Warburger Straße 100, 33098 Paderborn, Germany. Email: \href{mailto:sofya.maslovskaya@upb.de}{sofya.maslovskaya@upb.de}}\ \thanks{The work of this author has been supported by Deutsche Forschungsgemeinschaft (DFG), Grant No. OB 368/5-1, AOBJ: 692093}
\\
Sina Ober-Bl\"obaum\footnote{  Universität Paderborn (UPB), Numerical Mathematics and Control (NMC), Warburger Straße 100, 33098 Paderborn, Germany. Email: \href{mailto:sinaober@math.uni-paderborn.de}{sinaober@math.uni-paderborn.de}}\
\qquad
Rodrigo T.~Sato Mart{\'\i}n de Almagro\footnote{ Friedrich-Alexander-Universität Erlangen-Nürnberg (FAU), Institute of Applied Dynamics (LTD), Immerwahrstrasse 1, 91058 Erlangen, Germany. Email: \href{mailto:rodrigo.t.sato@fau.de}{rodrigo.t.sato@fau.de}}\
}

\begin{document}
\maketitle                   % Produces the title.

\begin{abstract}
It has been shown recently that optimal control problems with the dynamical constraint given by second order system admit a regular Lagrangian formulation. This implies that the optimality conditions can be obtained in a new form based on the variational approach. In this paper we extend the first order necessary optimality conditions obtained in  \cite{sato2025} to second order optimality conditions. This results in a complete characterization of the optimality conditions in a new Lagrangian form.
\end{abstract}
%\textcolor{black}{Maybe we can come up with a nicer title, e.g. Variational framework for ... instead of New framework} \sofya{Variational already exists and it is not fully variational, the proof mostly relies on the Tulczyjew's diffeo, maybe "New Lagrangian"?}
%\textcolor{black}{Fine for me: New Lagrangian framework...}

\section{Introduction}
{Optimal control plays an important role in modern engineering for the design of mechanical systems by optimising some objective function, which allows for e.g.\ the fuel- or time-efficient design of robots, satellites manoeuvres. }  Second order optimal control problems (OCP)s form a subclass where the control system is given by differential equations of second order. Such systems model in particular mechanical systems with external forces. In this case, physical laws lead to controlled Euler-Lagrange equations.

One of the main approaches to solve OCPs is based on the necessary conditions of optimality given by Pontryagin's maximum principle \cite{pontryagin1964,gerdts2023optimal}. Since the resulting equations are usually difficult to solve analytically, {numerical methods have to be} employed. This kind of numerical approaches to solve OCPs is called the indirect approach. Alternatively, the direct approach \cite{formalskii10,betts1999direct,gerdts2023optimal} is based on the dicretisation of the OCP in order to approximate it by a finite dimensional optimisation problem, and then to solve the corresponding necessary optimality conditions. A new framework for the class of second order OCPs  was recently introduced in \cite{leyendecker2024new,kono2025,sato2025,kono2025b}, where the OCP was first reformulated as a regular Lagrangian problem, and then the first order necessary optimality conditions for this problem were analysed. The advantage of the new framework is the possibility to use {theoretical tools from} 
%the methodology of 
discrete mechanics \cite{marsden2001discrete} 
%with the corresponding  
for the analysis of the discretised problem. However, the first order conditions do not allow to distinguish the minimising solutions among the extremals, this can be done by considering the second order conditions. This is why, this work has as a goal to extend the results from \cite{sato2025} by defining the second order optimality conditions in the new Lagrangian setting. %{Michael, I reformulated a bit the paragraph, to put more accent on the analytical solution of OCPs in the beginning (this paper is not on numerical methods) and also included our previous works.}{sounds fine that way. (originally i wanted to add our new stuff after talking about the standard ways, so would have added it after the conjugate point iscussion. But how it is now should be better.)}

%{Optimal control problems (OCP)s are usually solved by determining solutions of necessary optimality conditions, either in the form of indirect approaches which introduce the Pontryagin Hamilton function and utilize the Pontryagin's maximum principle \cite{pontryagin1964,gerdts2023optimal} or direct approaches \cite{formalskii10,betts1999direct,gerdts2023optimal} that discretize the OCP  and find solutions of the necessary optimality conditions of the subsequent optimization problem.}

{%While the necessary optimality conditions permit to find solution curves that fulfil the necessary conditions and thus are extremal, they are not guaranteed to be (locally) optimal. 
To find solutions that are guaranteed to be (locally) optimal, further sufficiency conditions need to be derived. Such sufficient conditions for local optima are given by the so-called  (strong) Legendre conditions \cite{agrachev2013control} which give second-order inequality conditions on the second control derivative of the Pontryagin's Hamilton function at the optimum. If the extremal solutions satisfy the strong Legendre condition, then there exists a time interval $[0,\epsilon]$ for which the solution curve is optimal.} 
The time when the extremal curves lose their optimality can be analysed in terms of the conjugate times. The theory of conjugate times was developed both in context of calculus of variations \cite{gelfand1963calculus} and in the context of optimal control \cite{Bonnard2007}. In both cases, the conjugate times are characterised by non trivial solutions of variational equations. Computation of conjugate points in  practice is usually done by simultaneous resolution of the extremal and variational equations in the indirect numerical approach. 
In this work we first complete the first order necessary optimality conditions obtained in \cite{sato2025} with the second order necessary conditions. In {particular}, we provide the {necessary and sufficient} optimality conditions equivalent to the classical Legendre conditions. Finally, we derive the equations on conjugate points in the new Lagrangian setting. These results lay the foundation for similar developments in the discrete setting.  %  obtained provide a general description of necessary and sufficient optimality condition for second order OCPs in the new Lagrangian framework. The new formulation in the setting of 

The paper is organised as follows. In Sec.~2 we recall the standard optimality conditions in optimal control. In Sec.~3 we state the first order optimality conditions from \cite{sato2025} in the new Lagrangian framework. We derive the additional second order necessary and sufficient optimality conditions in Sec.~4. In Sec.~5 the results are illustrated by an example of linear-quadratic OCP and we finish the paper with conclusions in Sec.~6.
%Optimal control problems of second order play especially important role in engeneering due to its connection to mechanical systems. 

%Optimal control problems are usually solved by applying necessary optimality conditions in form of Pontryagin's maximum principle. In the series of works \cite{leyendecker2024new,kono2025,sato2025,kono2025b} the Authors have shown that in case of second order optimal control problems, it is possible to use the Lagrangian formulation. It is especially helpful in numerical methods, because it allows to use methods of discrete mechanics \cite{marsden2001discrete}.

%In this work we give a general description of necessary optimality condition in the new Lagrangian framework.

\section{Standard optimality conditions for second order problems}
Following \cite{sato2025}, let us consider a state-control vector bundle $(\mathcal{E}, \pi^\mathcal{E}, \mathcal{M})$ where the base space is given by the state space $\mathcal{M} = T\mathcal{Q}$ with $\mathcal{Q}$ a smooth configuration manifold, and the control space by a vector space $\mathcal{N}$. We want to solve an OCP defined in local coordinates on $\mathcal{E}$ in the form
\begin{equation} \label{eq:OCP}
\begin{aligned}
	\min_{(q,u)} \quad &&      \mathcal{J}[q,u] &= \int_0^{T} C(q(t),\dot{q}(t),u(t))~dt\\
	\text{subject to} \quad &&        \ddot{q}(t) &= f(q(t), \dot{q}(t), u(t)),\\
	&& q(0) &= q_0, \quad \dot q(0) = v_0, \\
        && q(T) &= q_T, \quad \dot q(T) = v_T,
\end{aligned}
\end{equation}
where $q \in Q$ is a configuration, $u$ denotes the control, $q_0, v_0$ are initial configuration and velocity and $q_T, v_T$ are final configuration and velocity, cost function $C$ and vector field $f$ are assumed to be of class $C^2$. The necessary conditions of optimality for such optimal control problems are given by Pontryagin's Maximum Principle (PMP) \cite{pontryagin1962}. To apply PMP, \eqref{eq:OCP} should be first transformed into first order form by introducing a new variable $v$ for $\dot q$. This leads to the following OCP
\begin{equation} \label{eq:OCP.1ord}
	\begin{aligned}
		\min_{(q,v,u)} \quad &&      \mathcal{J}[q,v,u] &= \int_0^{T} C(q(t),v(t),u(t))~dt\\
		\text{subject to} \quad &&        \dot{q}(t) &= v(t),\\
        &&        \dot{v}(t) &= f(q(t), v(t), u(t)),\\
		&& q(0) &= q_0, \quad v(0) = v_0, \\
        && q(T) &= q_T, \quad v(T) = v_T.
	\end{aligned}
\end{equation}
  In particular, in the setting of problem \eqref{eq:OCP} with no constraint on the control, the weak PMP applies \cite{bonnard2003singular,trelat2024control}. The {Pontryagin's} Hamiltonian associated with \eqref{eq:OCP} is given by 
$$ \mathcal{H}_{\lambda^0}(x,\lambda, u) = {\langle (\lambda_q,\lambda_v), (v, f(q, v, u)) \rangle} +  \lambda^0 C(q,v,u),$$
{where $\lambda^0 \in \mathbb{R}$ and $\lambda = (\lambda_q,\lambda_v) \in T^*\mathcal{M}$ are called abnormal multiplier and adjoint (also costate), respectively}.
\begin{theorem}[Weak Pontryagin's Maximum Principle]
	Assume that $u(\cdot) \in L^\infty([0,T], \mathcal{N})$. Let the control $u(\cdot)$ and the corresponding state $x(\cdot)= (q(\cdot), v(\cdot))$ be optimal for \eqref{eq:OCP.1ord}. Then, there exits an absolutely continuous adjoint $\lambda(\cdot) = (\lambda_q(\cdot), \lambda_v(\cdot))$ and a constant $\lambda^0 \leq 0$ such that $(\lambda(\cdot), \lambda^0) \neq 0$ for any $t \in [0,T]$ and $(x,\lambda,\lambda^0, u)$ satisfy the following system of equations
    \begin{subequations}
    \label{eq:PMP}
        \begin{align}
		\dot q(t) &= \hphantom{-}\frac{\partial \mathcal{H}_{\lambda^0}}{\partial \lambda_q}(q(t),v(t),\lambda_q(t),\lambda_v(t), u(t)), & {\dot v(t)} &= \hphantom{-}\frac{\partial \mathcal{H}_{\lambda^0}}{\partial \lambda_v}(q(t),v(t),\lambda_q(t),\lambda_v(t), u(t)), \\
		\dot \lambda_q(t) &= -\frac{\partial \mathcal{H}_{\lambda^0}}{\partial q}(q(t),v(t),\lambda_q(t),\lambda_v(t), u(t)), & {\dot\lambda_v(t)} &= -\frac{\partial \mathcal{H}_{\lambda^0}}{\partial v}(q(t),v(t),\lambda_q(t),\lambda_v(t), u(t)), \\
		0 &= \hphantom{-}\frac{\partial \mathcal{H}_{\lambda^0}}{\partial u}(q(t),v(t),\lambda_q(t),\lambda_v(t), u(t)). && \label{eq:PMP_optimality}
	\end{align}
    \end{subequations}
and boundary conditions
\begin{subequations} \label{eq:PMP.boundary}
	\begin{alignat}{2}
		q(0) &= q_0, \quad & v(0) &= v_0, \\
		q(T) &= q_T, \quad & v(T) &= v_T.
        %\lambda_q(T) &= \lambda^0 \frac{\partial \phi}{\partial q}(q(T), v(T)), \quad \lambda_v = \lambda^0 \frac{\partial \phi}{\partial v}(q(T), v(T)).
	\end{alignat}
\end{subequations} 
\end{theorem}
{Solutions $(x, \lambda, \lambda^0,u)$ of \eqref{eq:PMP}-\eqref{eq:PMP.boundary} are candidates to be minimisers for \eqref{eq:OCP.1ord}. Given a candidate pair $(x,u)$, the dimension $m \geq 1$ of the vector space spanned by $(\lambda(T),\lambda^0)$ for all the admissible pairs $(\lambda(\cdot),\lambda^0)$ satisfying \eqref{eq:PMP} is said to be the corank of the minimiser. In this work it is assumed that we are always in the corank $m = 1$ case.}

The previous equations are invariant with respect to a rescaling of $(\lambda_q, \lambda_v, \lambda^0)$, therefore, it is classical to fix the value of $\lambda^0$ and distinguish two cases, namely, the normal case when $\lambda^0 = -1$, and the abnormal case when $\lambda^0 = 0$. In our consideration, we will restrict our attention to the normal case.

{If we define the augmented objective functional associated to \eqref{eq:OCP.1ord},
\begin{equation}
    \widetilde{\mathcal{J}}[x,\lambda,u] = \int_{0}^T \left[C(q(t),v(t),u(t) + \langle (\lambda_q(t),\lambda_v(t)), (\dot{q}(t) - v(t), \dot{v}(t) - f(q(t), v(t), u(t))) \rangle \right] \,dt, \label{eq:augmented.1ord}
\end{equation}
it should be noted that Equations \eqref{eq:PMP} imply the vanishing of the first variation of \eqref{eq:augmented.1ord}. That is, for all $(\delta x(\cdot),\delta \lambda(\cdot), \delta u(\cdot))$ such that $(\delta x(\cdot), \delta \lambda(\cdot)) \in T_{(x(\cdot),\lambda(\cdot))}T^*\mathcal{M}$ and $(\delta x(\cdot),\delta u(\cdot)) \in T_{(x(\cdot),u(\cdot))}\mathcal{E}$, we have that
\begin{equation*}
    d\widetilde{\mathcal{J}}_{[x,\lambda,u]}[\delta x,\delta \lambda,\delta u] = 0.
\end{equation*}}

The second order necessary optimality conditions allow us to distinguish minimising and maximising curves. The sufficient optimality conditions %have a local nature and
ensure optimality of a given curve. %for {sufficiently short time.}
In the case of \eqref{eq:OCP.1ord}, the {local} {necessary and sufficient} second order conditions are given by the Legendre conditions and can be summarised as follows \cite{Bonnard2007}.
\begin{theorem}[Legendre conditions \cite{Bonnard2007}] \label{th.Legendre} $\left. \right.$ 
\begin{itemize}
	\item If the control $u(\cdot)$ and the corresponding state $x(\cdot)$ are optimal for \eqref{eq:OCP.1ord}, then  $(x,\lambda,\lambda^0, u)$ satisfy the Legendre condition 
	\begin{equation} \label{eq:Legendre}
		\frac{\partial^2 \mathcal{H}_{{\lambda_0}}}{\partial u^2}(x(t),\lambda(t), u(t)) \leq 0, \qquad {\forall t\in[0,T].}
	\end{equation} 
    
	\item If $(x,\lambda,\lambda^0, u)$ satisfy the strong Legendre condition 
	\begin{equation} \label{eq:Legendre.strong}
		\frac{\partial^2 \mathcal{H}_{{\lambda_0}}}{\partial u^2}(x(t),\lambda(t), u(t)) <0, \qquad {\forall t\in[0,T],}
	\end{equation}  
then, there exists {$\varepsilon > 0$ such that $(x,u)$ is optimal throughout the time interval $[0,\varepsilon]$}. %time interval $[0,\varepsilon]$ on which $(x,u)$ is locally optimal.
\end{itemize}

\end{theorem}

{One can refine this result further introducing the notion of conjugate times. First, we can consider the second variation of \eqref{eq:augmented.1ord}, $d^2\tilde{\mathcal{J}}$. In the corank 1 case, when restricted to candidate minimisers, and fixing $\delta \lambda = 0$, this provides a unique quadratic form $Q = \left.d^2\tilde{\mathcal{J}}\right\vert_{\ker d\tilde{\mathcal{J}}}$ on the space of variations $T_{(x(\cdot),u(\cdot))}\mathcal{E}$, see \cite{Bonnard2007,agrachev2013control}. Theorem \ref{th.Legendre} guarantees the positive definiteness of $Q$ for sufficiently small intervals. The times at which the signature of $Q$ changes are called conjugate times and under the strong Legendre condition these can be shown to be isolated \cite{agrachev2013control}. The first of these conjugate times along a candidate minimiser marks the time where $Q$ stops being definite and becomes indefinite. This leads to the following
\begin{theorem}[Second order necessary and sufficient conditions for optimality \cite{agrachev2013control,Bonnard2007}]
    \label{th:PMP.second.order.refined}
    Let $(x(\cdot),\lambda(\cdot),u(\cdot))$ be a normal extremal of corank 1 of the optimal control problem \eqref{eq:OCP.1ord} on the interval $[0,T]$. Moreover, assume that the strong Legendre condition is satisfied along this extremal. Finally, let $t_c$ denote the first conjugate time along this extremal. Then, if $t_c > T$, the extremal is locally optimal\footnote{{The version cited assumes that the strong Pontryagin's maximum principle is satisfied and thus the extremal is locally strongly optimal (i.e. in the $C^0$ topology). However, in the case at hand, the extremal is only locally weakly optimal (i.e. in the $C^1$ topology).}}.
    Moreover, if the extremal is corank 1 on every subinterval of $[0,T]$, then, the extremal is not optimal if $t_c < T$.
\end{theorem}
In our case, the study of conjugate times can be carried out by examining the solutions of the variational equations of \eqref{eq:PMP}. In particular, let us assume an arbitrary differentiable family of solutions $(x_{\epsilon
}(\cdot),\lambda_{\epsilon
}(\cdot),u_{\epsilon
}(\cdot))$ of \eqref{eq:PMP}, with $\epsilon \in I$, and $I$ an open {interval} containing $0$. Further, assume that for $\epsilon = 0$, the corresponding element of the family $(x_{0}(\cdot),\lambda_{0}(\cdot),u_{0}(\cdot)) = (x(\cdot),\lambda(\cdot),u(\cdot))$ satisfies \eqref{eq:PMP.boundary}. Finally, let us denote $\delta x(\cdot) = \left.\frac{d}{d \epsilon}\right\vert_{\epsilon = 0} x(\cdot)$ and similarly for the rest of the variables. One can then insert this family into \eqref{eq:PMP}, differentiate with respect to $\epsilon$ and evaluate at $0$ to obtain the variational equations
    \begin{subequations}
    \label{eq:PMP.variational}
        \begin{align}
		\delta \dot{x} &= \hphantom{-}\frac{\partial^2 \mathcal{H}_{\lambda^0}}{\partial x \partial \lambda} \delta x + \frac{\partial^2 \mathcal{H}_{\lambda^0}}{\partial \lambda \partial \lambda} \delta \lambda + \frac{\partial^2 \mathcal{H}_{\lambda^0}}{\partial u \partial \lambda} \delta u, \label{eq:PMP.variational.x}\\
		\delta \dot{\lambda} &= -\frac{\partial^2 \mathcal{H}_{\lambda^0}}{\partial x \partial x} \delta x - \frac{\partial^2 \mathcal{H}_{\lambda^0}}{\partial \lambda \partial x} \delta \lambda - \frac{\partial^2 \mathcal{H}_{\lambda^0}}{\partial u \partial x} \delta u,\label{eq:PMP.variational.lambda}\\
		0 &= \hphantom{-}\frac{\partial^2 \mathcal{H}_{\lambda^0}}{\partial x \partial u} \delta x + \frac{\partial^2 \mathcal{H}_{\lambda^0}}{\partial \lambda \partial u} \delta \lambda + \frac{\partial^2 \mathcal{H}_{\lambda^0}}{\partial u \partial u} \delta u. \label{eq:PMP.variational.optimality}
	\end{align}
    \end{subequations}
This is a linear system of equations for the fields of variations along a given candidate. Those fields satisfying \eqref{eq:PMP.variational} are called Jacobi fields. Conjugate times can be found as times $t_c$ where there exist non-trivial solutions of the variational equations satisfying the boundary conditions
\begin{equation}
    \delta x(0) = 0, \quad \delta x(t_c) = 0.\label{eq:PMP.variational.BC}
\end{equation}}
 Under the strong Legendre conditions, the implicit function theorem can be applied to the optimality condition {\eqref{eq:PMP_optimality}} %$\frac{\partial \mathcal{H}_{\lambda^0}}{\partial u}(x,\lambda, u) = 0$
 in order to find an expression for the optimal control as a function of state and adjoint $u^* = u^*(x,\lambda)$. Substituting it in the Hamiltonian {gives us the \emph{so-called} reduced Hamiltonian,} $\mathcal{H}_{\lambda^0}^r(x,\lambda){:=}\mathcal{H}_{\lambda^0}(x,\lambda, u(x,\lambda))${, which} leads to a simplified formulation of \eqref{eq:PMP} as a standard Hamiltonian system on $T^*\mathcal{M}$ of the form 
\begin{equation} \label{eq:reduced.Ham}
 \dot x = \frac{\partial \mathcal{H}_{\lambda^0}^r}{\partial \lambda}(x(t), \lambda(t)), \qquad  \dot \lambda = - \frac{\partial \mathcal{H}_{\lambda^0}^r}{\partial x}(x(t), \lambda(t)).
\end{equation}
The solution $(x({\cdot}),\lambda({\cdot}))$ of the reduced Hamiltonian defined on $[0,T]$ that satisfies the boundary conditions \eqref{eq:PMP.boundary} is called an extremal. {It is not difficult to show that t}he linearisation of \eqref{eq:reduced.Ham} around an extremal 
\begin{equation} \label{eq:reduced.Ham.lin}
\begin{aligned}
	\delta \dot x &= \hphantom{-} \frac{\partial^2 \mathcal{H}_{\lambda^0}^r}{\partial x \partial \lambda}(x(t), \lambda(t)) \delta x + \frac{\partial^2 \mathcal{H}_{\lambda^0}^r}{\partial \lambda^2}(x(t), \lambda(t)) \delta \lambda,\\
    \delta \dot \lambda &= -\frac{\partial^2 \mathcal{H}_{\lambda^0}^r}{\partial x^2}(x(t), \lambda(t))\delta x -\frac{\partial^2 \mathcal{H}_{\lambda^0}^r}{\partial x \partial \lambda}(x(t), \lambda(t))\delta \lambda.
    \end{aligned}
\end{equation}
{coincides with \eqref{eq:PMP.variational.x}-\eqref{eq:PMP.variational.lambda}, over variations satisfying \eqref{eq:PMP.variational.optimality}. This provides a lower-dimensional system to determine conjugate times.}

{The relation}
between optimality and conjugate points permits us to establish algorithms to check optimality of an extremal, as it has been done in \cite{Bonnard2007}. {Their algorithm} is based on the integration of \eqref{eq:reduced.Ham.lin} from $n=\dim \mathcal{M}$ initial conditions $(\delta x_i(0), \delta  \lambda_i(0)) = (0, e_i)$, for $i = 1, \dots n$, where $e_i$ is the $i$th basis vector in $\R^n$. Then the conjugate time $t_c$ is found as the first time $t$ at which there holds $\det(\delta x_1(t), \cdots, \delta x_n(t)) = 0$.

\section{First order conditions in the new Lagrangian formulation}

The first order optimality conditions equivalent to \eqref{eq:PMP}-\eqref{eq:PMP.boundary} can be obtained without transformation to {first order form as in} \eqref{eq:OCP.1ord} using the variational approach under more restrictive regularity assumptions on $u$, namely, $u \in C([0,T], \mathcal{N})$. The variational approach is based on the consideration of the augmented Lagrangian
\begin{equation} \label{aug.Lagrangian}
	\mathcal{J}_{\mathrm{aug}}[q, \kappa ,u] = \int_0^{T} C(q(t),\dot{q}(t),u(t))+ \kappa^{\top}(t) \left(\ddot q(t) -  f(q(t), \dot{q}(t), u(t)) \right)~dt.
\end{equation}
Assume $q, \kappa$ are of class $C^2$ and $u$ of class $C$. Then imposing the vanishing of the first variation of $\mathcal{J}_{\mathrm{aug}}$ leads to the following equations
\begin{equation} \label{eq:aug.Lagrangian}
	\begin{aligned}
		\ddot q(t) &= {f(q(t), \dot{q}(t), u(t))}, \\
		\ddot \kappa(t) &= \frac{d}{dt} \left[\frac{\partial C}{\partial \dot q}(q(t),\dot{q}(t),u(t)) - \kappa(t)^{\top} \frac{\partial f}{\partial \dot q}(q(t),\dot{q}(t),u(t))\right] - \frac{\partial C}{\partial q} (q(t),\dot{q}(t),u(t)) + \kappa(t)^{\top} \frac{\partial f}{\partial q}(q(t),\dot{q}(t),u(t)),\\
		0 &= \frac{\partial C}{\partial u}(q(t),\dot{q}(t),u(t)) - \kappa(t)^{\top} \frac{\partial f}{\partial u}(q(t),\dot{q}(t),u(t)),
	\end{aligned}
\end{equation} 
and boundary conditions
\begin{equation} \label{eq:aug.Lagrangian.boundary}
	\begin{aligned}
		q(0) &= q_0, \quad v(0) = v_0,  \\
		q(T) &= q_T, \quad v(T) = v_T.
        %\kappa(T) &=\frac{\partial \phi}{\partial q}(q(T),\dot{q}(T)), \ \dot{\kappa}(T)^{\top} = \frac{\partial \phi}{\partial q}(q(T),\dot{q}(T)) +\frac{\partial C}{\partial \dot q}(q(T),\dot{q}(T),u(T)) + \frac{\partial \phi}{\partial \dot q}(q(T),\dot{q}(T)) \, \frac{\partial f}{\partial \dot q}(q(T),\dot{q}(T),u(T)).
	\end{aligned}
\end{equation} 
As shown in \cite{sato2025}, \eqref{eq:aug.Lagrangian}-\eqref{eq:aug.Lagrangian.boundary} are equivalent to \eqref{eq:PMP}-\eqref{eq:PMP.boundary} under {the}  identification $\lambda_v = \kappa$ and $\lambda_q = \frac{\partial C}{\partial \dot q}(q,\dot{q},u) - \langle \kappa, \frac{\partial f}{\partial \dot q}(q,\dot{q},u) - \dot{\kappa}\rangle$.

The new approach to first order optimality conditions introduced in \cite{leyendecker2024new,kono2025,sato2025,kono2025b} is based on the new Lagrangian formulation obtained from \eqref{aug.Lagrangian} by integration by parts
\begin{multline} \label{new.Lagrangian.cost}
	\widetilde{ \mathcal{J}}_{\mathrm{aug}}[q, \kappa ,u] =  \kappa(T)^{\top} \dot{q}(T) - \kappa(0)^{\top} \dot{q}(0)  + \int_{0}^T \left[ C(q(t),\dot{q}(t),u(t)) - \dot{\kappa}(t)^{\top} \dot{q}(t) - \kappa(t)^{\top} f(q(t),\dot{q}(t),u(t))\right] dt.
\end{multline}
The new control Lagrangian $\tilde{\mathcal{L}}^\mathcal{E}$ is defined on $TT^*\mathcal{Q} \oplus_{T\mathcal{Q}}^{\alpha} \mathcal{E}$ with values in $\mathbb{R}$, see \cite{sato2025} for more details on the geometrical setting. Notation $TT^*\mathcal{Q} \oplus_{T\mathcal{Q}}^\alpha \mathcal{E}$ stands for the $\alpha_\mathcal{Q}$-twisted sum of $TT^*\mathcal{Q}$ and $\mathcal{E}$ defined as follows
$$TT^*\mathcal{Q} \oplus_{T\mathcal{Q}}^\alpha \mathcal{E} = \left\lbrace (V,U) \in TT^*\mathcal{Q} \times \mathcal{E} \;\vert\; \pi_{T\mathcal{Q}} \circ \alpha_{\mathcal{Q}} (V) = \pi^{\mathcal{E}}(U) \right\rbrace,$$ 
with $\alpha_{\mathcal{Q}}$ Tulczyjew’s isomorphism \cite{TulczyjewLag} between the double bundles $T^*T\mathcal{Q}$ and $TT^*\mathcal{Q}$. In local coordinates $(q,\kappa,v_q,v_{\kappa})$ on $T^*T\mathcal{Q}$ it is defined by
\begin{equation*}
	\alpha_{\mathcal{Q}}(q,\kappa,v_q,v_{\kappa}) = (q,v_q,v_{\kappa},\kappa) \in T^*T\mathcal{Q}.
%	\beta_{\mathcal{Q}}(q,\kappa,v_q,v_{\kappa}) &= (q,\kappa,-v_{\kappa},v_q) \in T^*T^*\mathcal{Q},
\end{equation*}
In local coordinates $(q,\kappa,v_q,v_{\kappa},u)$ on $TT^*\mathcal{Q} \oplus_{T\mathcal{Q}}^{\alpha} \mathcal{E}$, the new Lagrangian is given by 
\begin{equation} \label{new.Lagrangian}
	\tilde{\mathcal{L}}^\mathcal{E}(q,\kappa,v_q,v_{\kappa},u) = v_{\kappa}^{\top} v_q + \kappa^{\top} f(q,v_q,u) - C(q,v_q,u)\,.
\end{equation}
For {the sake of notational economy}, we introduce the notation $y = (q, \kappa)$. As it was shown in \cite{sato2025}, the matrix $\frac{\partial^2 \tilde{\mathcal{L}}^\mathcal{E}}{\partial \dot y^2}(y, \dot y, u)$ is non degenerate, which implies that the new control Lagrangian is hyperregular and the first order conditions of optimality take the form of regular Euler-Lagrange equations involving the control variable as an additional parameter.

The relation between the new Lagrangian formulation and the Hamiltonian formulation from PMP was established in \cite{sato2025} based on an extension of Tulczyjew’s isomorphism. 
%The classical Tulczyjew’s isomorphism \cite{TulczyjewLag} is defined as an isomorphism between the double bundles $T^*T\mathcal{Q}$ and $TT^*\mathcal{Q}$. Let $(q,\kappa,v_q,v_{\kappa})$ be local coordinates on $T^*T\mathcal{Q}$. Then we have the following relation
%\begin{equation*}
%	\alpha_{\mathcal{Q}}(q,\kappa,v_q,v_{\kappa}) = (q,v_q,v_{\kappa},\kappa) \in T^*T\mathcal{Q}.
%	\beta_{\mathcal{Q}}(q,\kappa,v_q,v_{\kappa}) &= (q,\kappa,-v_{\kappa},v_q) \in T^*T^*\mathcal{Q},
%\end{equation*}
%This formulation was extended in \cite{leyendecker2024new,sato2025} to double control bundles. 
In adapted local coordinates, it is defined as follows%Tulczyjew's diffeomorphism has the following form
\begin{equation*}
	\alpha^\mathcal{E}_{\mathcal{Q}}(q,\kappa,v_q,v_{\kappa}, u) = (q,v_q,v_{\kappa},\kappa, u) \in T^*T\mathcal{Q}\oplus_{T\mathcal{Q}} \mathcal{E},
%	\beta^\mathcal{E}_{\mathcal{Q}}(q,\kappa,v_q,v_{\kappa}, u) &= (q,\kappa,-v_{\kappa},v_q, u) \in T^*T^*\mathcal{Q}\oplus_{T\mathcal{Q}} \mathcal{E},
\end{equation*}
with $(q,\kappa,v_q,v_{\kappa}) \in TT^*\mathcal{Q}$. It induces the following relation
\begin{align} \label{eq:Tulc}
	\tilde{\mathcal{L}}^\mathcal{E} &= \mathcal{H}_{\lambda^0} \circ \alpha_{\mathcal{Q}}^{\mathcal{E}}. % \qquad \text{or equivalently} \qquad \tilde{\mathcal{L}}^\mathcal{E}\circ (\alpha_{\mathcal{Q}}^{\mathcal{E}})^{-1} = \mathcal{H}_{\lambda^0}.
\end{align}
As a result, the extended Tulczyjew's diffeomorphism makes the relation between the optimality conditions given by the weak PMP and obtained from the variational approach in the new Lagrangian setting {explicit}. In the next section we will further describe the implications on the second order necessary and sufficient optimality conditions.
\section{Second order conditions in the new Lagrangian setting}

Assume that $(q^*, \kappa^* ,u^*)$ satisfy the first order optimality conditions \eqref{eq:aug.Lagrangian} and $(q^* ,u^*)$ minimise \eqref{eq:OCP}. Let $(q,u)$ satisfy $\ddot q = f(q,\dot q, u)$, then the following inequality holds
\begin{equation} \label{eq:frozen.adj}
	\widetilde{ \mathcal{J}}_{\mathrm{aug}}[q, \kappa^* ,u] - \widetilde{ \mathcal{J}}_{\mathrm{aug}}[q^*, \kappa^* ,u^*] = \mathcal{J}[q, u] - \mathcal{J}[q^*, u^*] \geq 0.
\end{equation}
%for all $(q,u)$ satisfying $\ddot q = f(q,\dot q, u)$. 
{In the following, w}e denote $\widetilde{ \mathcal{J}}_{\mathrm{aug}}^* [q ,u] = \widetilde{ \mathcal{J}}_{\mathrm{aug}}[q, \kappa^* ,u]$. {From}
%It follows from
\eqref{eq:frozen.adj} {it follows}, that $(q^*, u^*)$ minimise $\widetilde{ \mathcal{J}}_{\mathrm{aug}}^* [q ,u]$ under the condition $\ddot q = f(q,\dot q, u)$, and thus, its second variation must be non-negative. Using the notation $x = (q,\dot q)$, we have
$$ 
0 \leq \delta^2\widetilde{ \mathcal{J}}_{\mathrm{aug}}^* =  - \int_0^T (\delta x^\top, \delta u^\top) \left(
\begin{array}{cc}
	\frac{\partial^2 \tilde{\mathcal{L}}^\mathcal{E} }{\partial x^2} & \frac{\partial^2 \tilde{\mathcal{L}}^\mathcal{E} }{\partial x \partial u} \\[0.1cm]
	\frac{\partial^2 \tilde{\mathcal{L}}^\mathcal{E} }{\partial x \partial u} & \frac{\partial^2 \tilde{\mathcal{L}}^\mathcal{E} }{\partial u^2} 
\end{array}
\right) \left(
\begin{array}{cc}
	\delta x \\
\delta u
\end{array}
\right) dt,
$$
where $\delta x$ and $\delta u$ are related by $\delta \dot x  = \frac{\partial f}{\partial x}(x^*,u^*)\delta x + \frac{\partial f}{\partial u}(x^*,u^*)\delta u$. Using the classical argument of the dominant term, see \cite[Section 3.4.3]{Liberzon2012}, we conclude that $-\frac{\partial^2 \tilde{\mathcal{L}}^\mathcal{E} }{\partial u^2}({y^*(t), \dot y^*(t)},u^*(t)) \geq 0$ for all $t \in [0,T]$. We also provide a proof for completeness. 
\begin{theorem}
    Let {$(y^*(\cdot), u^*(\cdot)) = (q^*(\cdot), \kappa^*(\cdot), u^*(\cdot))$} be an optimal solution of \eqref{eq:OCP} {with} $\kappa^*$ the corresponding adjoint from \eqref{aug.Lagrangian}. Then the following condition is satisfied
    $$ \frac{\partial^2 \tilde{\mathcal{L}}^\mathcal{E} }{\partial u^2}({y^*(t), \dot y^*(t)}, u^*(t)) \leq 0, \qquad t \in [0,T].$$
\end{theorem}
\begin{proof}
    The proof is by contradiction and follows the classical ideas in the calculus of variations for the proof of the Legendre necessary condition of optimality \cite{gelfand1963calculus}. Assume that there is $t_0 \in (0,T)$ such that $\frac{\partial^2 \tilde{\mathcal{L}}^\mathcal{E} }{\partial u^2}({y^*(t_0), \dot y^*(t_0)},u^*(t_0)) > 0$. Then there exists an interval $(t_0-h,t_0+h)$ such that $\frac{\partial^2 \tilde{\mathcal{L}}^\mathcal{E} }{\partial u^2}({y^*(t), \dot y^*(t)},u^*(t)) > \varepsilon$ for some $\varepsilon > 0$. We will choose variations $\delta q, \delta u$ such that $0 > \delta^2\widetilde{ \mathcal{J}}_{\mathrm{aug}}^*$, which would contradict optimality of $(q^*(\cdot), u^*(\cdot))$. The relation between $\delta q$ and $\delta u$ is given by
    $$\delta \ddot q = A_1(t) \delta q + A_2(t) \delta \dot q + B(t) \delta u,$$
    where matrices $A_1,A_2,B$ are defined as follows
    $$\quad A_1(t)= \frac{\partial f}{\partial q}(q^*(t), \dot q^*(t), u^*(t)), \ A_2(t)= \frac{\partial f}{\partial \dot q}(q^*(t), \dot q^*(t), u^*(t)), \ B(t)= \frac{\partial f}{\partial u}(q^*(t), \dot q^*(t), u^*(t)).$$
    We assume that $\mathrm{rank} B(t) = \dim \mathcal{N}$. Let us denote $\tilde B = (B^\top B)^{-1}B^\top$. Then we can express $\delta u$ as a function of $\delta q, \delta \dot q, \delta \ddot q$ as follows
    \begin{equation} \label{eq:deltau}
    \delta u = \tilde B(t) \delta \ddot q - \tilde B(t) A_1(t) \delta q - \tilde B(t) A_2(t) \delta \dot q.    
    \end{equation}
    We define variations $\delta q_i$ for $i = 1, \dots, n= \dim \mathcal{Q}$ of the form
    $$
    \delta q_i = \begin{cases}
        0, \quad 0 \leq t < t_0- h, \\
        \frac{6}{h^5}(t - (t_0 - h))^5 - \frac{15}{h^4}(t-(t_0-h))^4 + \frac{10}{h^3}(t-(t_0-h))^3, \quad  t_0- h\leq t \leq t_0, \\
        \frac{6}{h^5}((t_0 + h)-t)^5 - \frac{15}{h^4}((t_0 + h)-t)^4 + \frac{10}{h^3}((t_0 + h)-t)^3, \quad  t_0 \leq t \leq t_0+h, \\
        0, \quad t_0+h <  t \leq T.
    \end{cases}
    $$
    It follows from the construction that $\delta q = (\delta q_1, \dots, \delta q_n) \in C^2[0,T]$ and $\delta q $ satisfies the boundary conditions $\delta q(0) = \delta q(T) = 0$. In addition, we have the following
     $$
    \max_{t\in [0,T]} \|\delta q(t) \|_2  = \sqrt{n},  \  \max_{t\in [0,T]} \|\delta \dot q(t)\|_2  = \sqrt{n}\frac{15}{8h}, \  \max_{t\in [0,T]} \|\delta \ddot q_i(t)\|_2  = \sqrt{n}\frac{10 \sqrt{3}}{3h^2}.
    $$
    Together with \eqref{eq:deltau}, {this} implies
    $$
    \max_{t\in [0,T]} \|\delta u(t) \|_2 \leq M \sqrt{n} \left( 1 + \frac{15}{8h} + \frac{10 \sqrt{3}}{3h^2} \right), \quad M = \max\{ \max_{[0,T]}\tilde B(t), \max_{[0,T]}\tilde B(t) A_1(t), \max_{[0,T]}\tilde B(t)A_2(t)\}.
    $$
    The second variation $\delta^2\widetilde{ \mathcal{J}}_{\mathrm{aug}}^*$ corresponding to the chosen $\delta q$ has the following form
    \begin{multline}
         \delta^2\widetilde{ \mathcal{J}}_{\mathrm{aug}}^* =  -\int_{t_0-h}^{t_0+h} \delta q^\top(t)\frac{\partial^2 \tilde{\mathcal{L}}^\mathcal{E}_* }{\partial q^2} \delta q(t) + 2 \delta q^\top(t)\frac{\partial^2 \tilde{\mathcal{L}}^\mathcal{E}_* }{\partial q \partial \dot q} \delta \dot q(t) + \delta  \dot q^\top(t)\frac{\partial^2 \tilde{\mathcal{L}}^\mathcal{E}_* }{\partial \dot q^2}\delta \dot q(t) \\ + 2 \delta  q^\top(t)\frac{\partial^2 \tilde{\mathcal{L}}^\mathcal{E}_* }{\partial q \partial u} \delta u(t)  + 2 \delta  \dot q^\top(t)\frac{\partial^2 \tilde{\mathcal{L}}^\mathcal{E}_* }{\partial \dot q \partial u} \delta u(t) + \delta  u^\top(t)\frac{\partial^2 \tilde{\mathcal{L}}^\mathcal{E}_* }{\partial u^2} \delta u(t) dt, \end{multline}
         where we use the following notation for shortness $\tilde{\mathcal{L}}^\mathcal{E}_*(t) = \tilde{\mathcal{L}}^\mathcal{E}({y^*(t), \dot y^*(t)},u^*(t))$.
         We split $\delta^2\widetilde{ \mathcal{J}}_{\mathrm{aug}}^*$ into two terms $\delta^2\widetilde{ \mathcal{J}}_{\mathrm{aug}}^* = J_1(h) + J_2(h)$ defined as follows
         \begin{multline}
             J_1(h) = -\int_{t_0-h}^{t_0+h} \delta q^\top(t)\frac{\partial^2 \tilde{\mathcal{L}}^\mathcal{E}_* }{\partial q^2} \delta q(t) + 2 \delta q^\top(t)\frac{\partial^2 \tilde{\mathcal{L}}^\mathcal{E}_* }{\partial q \partial \dot q} \delta \dot q(t) + \delta  \dot q^\top(t)\frac{\partial^2 \tilde{\mathcal{L}}^\mathcal{E}_* }{\partial \dot q^2} \delta \dot q(t) \\ + 2 \delta  q^\top(t)\frac{\partial^2 \tilde{\mathcal{L}}^\mathcal{E}_* }{\partial q \partial u} \delta u(t)  + 2 \delta  \dot q^\top(t)\frac{\partial^2 \tilde{\mathcal{L}}^\mathcal{E}_* }{\partial \dot q \partial u} \delta u(t)dt,
         \end{multline}
         and 
         $$ J_2(h) = -\int_{t_0-h}^{t_0+h} \delta  u^\top(t)\frac{\partial^2 \tilde{\mathcal{L}}^\mathcal{E}_* }{\partial u^2} \delta u(t) dt.$$
         We have the following estimates on $J_1(h), J_2(h)$
         \begin{align}
             J_1(h) &\leq C_1h + C_2 +C_{-1}\frac{1}{h} + C_{-2} \frac{1}{h^2}, \\
             J_2(h) &\leq -\varepsilon \left(D_1h + D_2 +D_{-1}\frac{1}{h} + D_{-2} \frac{1}{h^2} + D_{-3} \frac{1}{h^3} \right)< 0,
         \end{align} 
         where the constants $C_i, D_i$ are positive and independent from $h$. This implies that there exists $0<h_0<1$ such that for all $0<h<h_0$ we have $| J_2(h)| > |J_1(h)|$, which implies $\delta^2\widetilde{ \mathcal{J}}_{\mathrm{aug}}^* < 0$. This contradicts the optimality of $(q^*(\cdot),u^*(\cdot))$. Therefore, we conclude that $\frac{\partial^2 \tilde{\mathcal{L}}^\mathcal{E} }{\partial u^2}({y^*(t), \dot y^*(t)},u^*(t)) \leq 0$ for any $ t \in [0,T]$.
\end{proof}

%{the equations are exactly as in \cite[Section 3.4.3]{Liberzon2012} , but a rigorous proof is probably required, as in \cite[Section 3.4.3]{Liberzon2012}, it is not formulated as a theorem. Notice that the same approach doesn't work when applied to the control independent Lagrangian. I tried to use it to get the equation on conjugate points, but the equation is wrong. At the same time $(y, \dot y) = (q, \kappa, \dot q, \dot \kappa)$ form a saddle point for the control independent Lagrangian, which follows from the indefiniteness of the non degenerate form $	\frac{\partial^2 \tilde{\mathcal{L}}}{\partial \dot y^2}(y, \dot y)$.}
% {Maybe naive, but is that saddle point related to the fact that the Tulczyjew transformation relates the costates with a '-' sign (so due to the partial integration)? And aren't you missing the partial integration terms here? } {The saddle point is due to the fact that our Lagrangian comes from adding the constraints. For example, in the finite dimensional minimization, for a Lagrangian formed from a cost function and constraints the pair $(x^*,\lambda^*)$, with $x^*$ optimal state and $\lambda^*$ the corresponding Lagrange multiplier, is a saddle point for the Lagrangian, so here it must be the same phenomenon behind. As for the partial integration terms, it is not present, because the variations are only in the direction of $q,\dot q$, so the second variation of $\kappa^\top \dot q$ vanishes.}
To obtain  {necessary and sufficient} second order conditions {in} the new Lagrangian setting, we {can simply} use Tulczyjew's diffeomorphisms {to translate those found earlier}. Using \eqref{eq:Tulc}, conditions in Theorem~\ref{th.Legendre} can be formulated in terms of the new Lagrangian, leading to the statement.
\begin{theorem} \label{th.Legendre.new} $\left. \right.$ 
	\begin{itemize}
		\item If the control $u(\cdot)$ and the corresponding state $q(\cdot)$ are optimal for \eqref{eq:OCP}, then  $({y, \dot{y}}, u)$ satisfy  the condition
		\begin{align} \label{eq:Legendre.new}
			\frac{\partial^2 \tilde{\mathcal{L}}^\mathcal{E}}{\partial u^2}({y(t), \dot y(t)}, u(t)) \leq 0 \qquad {\forall t\in [0,T].}
		\end{align} 
		\item If $({y, \dot{y}}, u)$ satisfy  the strong Legendre condition
		\begin{align} \label{eq:Legendre.strong.new}
			\frac{\partial^2 \tilde{\mathcal{L}}^\mathcal{E}}{\partial u^2}({y(t), \dot y(t)}, u(t)) <0, \qquad {\forall t\in [0,T],}
		\end{align}  
		then, there exists time interval $[0,\varepsilon]$ on which $(x,u)$ are locally optimal.
	\end{itemize}
\end{theorem}
The same approach can be used to obtain the equation \eqref{eq:reduced.Ham.lin} in terms of the new Lagrangian. Let us now assume that $\frac{\partial^2 \tilde{\mathcal{L}}^\mathcal{E} }{\partial u^2}({y^*(t), \dot y^*(t)},u^*(t)) < 0$ for all $t \in [0,T]$ and $u^* = u^*({y^*(t), \dot y^*(t)})$. We denote $\tilde{\mathcal{L}}(y,\dot y) = \tilde{\mathcal{L}}^\mathcal{E}(y,\dot y,u(y, \dot y))$. The optimal $y^*(\cdot)$ satisfies the Euler-Lagrange equation of the form
\begin{equation*} \frac{\partial \tilde{\mathcal{L}} }{\partial y} - \frac{d}{dt}\left(\frac{\partial \tilde{\mathcal{L}} }{\partial \dot y}\right) = 0.
\end{equation*}%\frac{\partial \tilde{\mathcal{L}} }{\partial y}  - \frac{\partial^2 \tilde{\mathcal{L}} }{\partial y \partial \dot y} \dot y - \frac{\partial^2 \tilde{\mathcal{L}} }{\partial \dot y^2} \ddot y = 0. $$
{The} linearised Euler-Lagrange equation along $y^*(\cdot)$ is defined as follows
\begin{equation}  \label{lin.EL}
	%0 &=  \frac{\partial^2 \tilde{\mathcal{L}} }{\partial y^2}\delta y +  \frac{\partial^2 \tilde{\mathcal{L}} }{\partial y \partial \dot y}\delta \dot y - \frac{\partial^2 \tilde{\mathcal{L}} }{\partial y \partial \dot y} \delta \dot y - \frac{\partial^3 \tilde{\mathcal{L}} }{\partial y^2 \partial \dot y}\dot y \delta y - \frac{\partial^3 \tilde{\mathcal{L}} }{\partial y \partial \dot y^2}\dot y \delta \dot  y - \frac{\partial^2 \tilde{\mathcal{L}} }{\partial \dot y^2} \delta \ddot y - \frac{\partial^3 \tilde{\mathcal{L}} }{\partial y \partial \dot y^2} \ddot y \delta y - \frac{\partial^3 \tilde{\mathcal{L}} }{\partial \dot y^3} \ddot y \delta \dot y \nonumber\\
   % & =  \frac{\partial^2 \tilde{\mathcal{L}} }{\partial y^2}\delta y +  \frac{\partial^2 \tilde{\mathcal{L}} }{\partial y \partial \dot y}\delta \dot y - \frac{d}{dt} \left( \frac{\partial^2 \tilde{\mathcal{L}} }{\partial y \partial \dot y} \delta  y +  \frac{\partial^2 \tilde{\mathcal{L}} }{\partial \dot y^2} \delta \dot y \right) \nonumber\\
     \left( \frac{\partial^2 \tilde{\mathcal{L}} }{\partial y^2} -  \frac{d}{dt} \frac{\partial^2 \tilde{\mathcal{L}} }{\partial y \partial \dot y}  \right) \delta  y - \frac{d}{dt} \left( \frac{\partial^2 \tilde{\mathcal{L}} }{\partial \dot y^2} \delta \dot y \right) = 0,
\end{equation}
which is a second-order equation in $\delta y$, known as Jacobi equation in the calculus of variations \cite{gelfand1963calculus}. {Similarly to how solutions of \eqref{eq:reduced.Ham.lin} together with \eqref{eq:PMP.variational.BC} give us the conjugate times, \eqref{lin.EL} together with \eqref{eq:PMP.variational.BC} give us conjugate times as well.}
% \begin{definition} \label{th.conjugate.points}
% 	Let $\delta y = (\delta q, \delta \kappa)$ be a nontrivial solution of the Jacobi equation satisfying $(\delta q(0), \delta \dot q(0)) = (0,0)$. The conjugate time is defined as $t_c>0$ such that $(\delta q(t_c), \delta \dot q(t_c)) = (0,0)$.
% \end{definition}
\begin{theorem}
	{The conjugate times obtained from \eqref{eq:reduced.Ham.lin} and \eqref{lin.EL}} coincide.
\end{theorem}
\begin{proof}
	{The} Hamiltonian system \eqref{eq:reduced.Ham} can be equivalently written in the form of second order differential equations in variables $(q, \lambda_v)$ identical to \eqref{eq:aug.Lagrangian} with $\lambda_v$ playing the role of $\kappa$, see \cite{sato2025} for details. As a result, \eqref{eq:reduced.Ham.lin} and \eqref{lin.EL} are equivalent and their solutions can be identified using $\delta \lambda_v = \delta \kappa$.
\end{proof}

{All of this permits one to apply Theorem \ref{th:PMP.second.order.refined} on the new control Lagrangian setting under the substitution $(x(\cdot),\lambda(\cdot),u(\cdot)) \mapsto (y(\cdot),u(\cdot))$.}

\section{Example}
Let us apply the derived theory to an example. We consider $\mathcal{Q} = \R^n$ and $\mathcal{M} = \R^{2n}$, the state-control bundle given by $\mathcal{E} = \R^{2n} \times \R^m$ and a linear-quadratic optimal control problem of the form
\begin{equation} 
		\begin{aligned}
				\mathcal{J}[q,u] &= \frac{1}{2} \int_0^T  q(t)^\top Q_1 q(t) + \dot q(t)^\top Q_2 \dot q(t) + u^\top(t) R u(t) dt\\
				\text{subject to}\ \ & \ddot q = A_1 q + A_2 \dot q + B u,\\
				& q(0) =q_0, \quad \dot q(0) = v_0, \\
                & q(T) =q_T, \quad \dot q(T) = v_T.
			\end{aligned}
		\label{LQ}
\end{equation}
{The} matrices $Q_1, Q_2 \in M_n(\R), \ R \in  M_m(\R)$ are assumed to be positive-definite and $A_1, A_2 \in  M_n(\R)$ and $B \in M_{n,m}(\R)$ are such that the following pair of matrices  
$$ \tilde{A} = \left(
\begin{array}{cc}
	0_{n,n} & Id_{n} \\
A_1 & A_2
\end{array}
\right), \qquad \tilde B = \left(
\begin{array}{cc}
	0_{n,m} \\
	B
\end{array}
\right), $$
with $0_{n,m} \in M_{n,m}$ {being a} matrix with zero entries, satisfies the Kalman condition, namely matrix $\left(\tilde B, \ \tilde A \tilde B, \cdots,  \tilde A^{2n-1} \tilde B \right)$ has rank $2n$. These assumptions ensure the controllablity of the system and the existence of an optimal solution \cite{trelat2024control}. Let us first write down the classical optimality conditions. The normal Hamiltonian has the form $\mathcal{H}_{\lambda^0} = \lambda_q^\top v  + \lambda_v^\top (A_1q + A_2 v + Bu) - \frac{1}{2} \left( q^\top Q_1 q + v^\top Q_2 v + u^\top R u \right)$ {and the} strong Legendre condition is satisfied {due to}  $\frac{\partial^2 \mathcal{H}_{\lambda^0}}{\partial u^2}(x(t),\lambda(t), u(t)) = -R <0$. The optimality condition $\frac{\partial \mathcal{H}_{\lambda^0}}{\partial u} (x(t),\lambda(t), u(t)) = 0$ allows {us} to express $u = R^{-1}B^\top \lambda_v$. The reduced Hamiltonian  is $\mathcal{H}^r_{-1} = \lambda_q^\top v  + \lambda_v^\top (A_1q + A_2 v)  - \frac{1}{2} \left( q^\top Q_1 q + v^\top Q_2 v + \lambda_v^\top B  R^{-1} B^\top \lambda_v \right)$ {resulting in a corresponding linear system}
%and the corresponding system is linear 
of the form
\begin{equation}
\begin{aligned} \label{eq:LQ.Ham.syst}
\dot q &= v, \qquad\dot v = A_1q + A_2 v + BR^{-1}B^\top \lambda_v, \\
\dot \lambda_q &= Q_1 q - A_1^\top \lambda_v, \qquad \dot \lambda_v = Q_2v-\lambda_q - A_2^\top \lambda_v.
\end{aligned}
\end{equation}
%Notice that the system is controllable, this is why we can use the variational equation to determine conjugate points. 
As the system is linear, linearisation around any extremal coincides with the system itself, namely
\begin{equation}
\begin{aligned} \label{eq:LQ.Ham.syst.var}
\delta \dot q &= \delta v, \qquad \delta \dot v = A_1\delta q + A_2 \delta v + BR^{-1}B^\top \delta \lambda_v, \\
\delta \dot \lambda_q &= Q_1 \delta q - A_1^\top \delta \lambda_v, \qquad \delta \dot \lambda_v = Q_2 \delta v-\delta \lambda_q - A_2^\top \delta \lambda_v.
\end{aligned}
\end{equation}
Let us now consider the new Lagrangian framework. The new control Lagrangian has the form $\tilde{\mathcal{L}}^\mathcal{E}(q,\kappa,v_q,v_{\kappa},u) = v_{\kappa}^{\top} v_q + \kappa^{\top} (A_1q + A_2 v + Bu) - \frac{1}{2} \left( q^\top Q_1 q + v^\top Q_2 v + u^\top R u \right)$. In this case, in accordance with the relationship between $\mathcal{H}_{\lambda^0}$ and {$\tilde{\mathcal{L}}^\mathcal{E}$}, the strong Legendre condition is satisfied $\frac{\partial^2 \tilde{\mathcal{L}}^\mathcal{E}}{\partial u^2}(x(t),\lambda(t), u(t)) = -R <0$. The new uncontrolled Lagrangian has the form $\tilde{\mathcal{L}}(q,\kappa,v_q,v_{\kappa}) = v_{\kappa}^{\top} v_q + \kappa^{\top} (A_1q + A_2 v) - \frac{1}{2} \left( q^\top Q_1 q + v^\top Q_2 v - \kappa^\top B  R^{-1} B^\top \kappa  \right)$ and the corresponding Euler-Lagrange equations read
\begin{equation}
\begin{aligned} \label{eq:LQ.EL.syst}
\ddot q &= A_1 q + A_2 \dot q + BR^{-1}B^\top \kappa, \\
\ddot \kappa &= \left(Q_2A_1-Q_1\right)q + Q_2A_2\dot q + \left( A_1^\top + Q_2BR^{-1}B^\top \right)\kappa - A_2^\top \dot \kappa.
\end{aligned}
\end{equation}
The corresponding Jacobi equation has the form
\begin{equation}
\begin{aligned} \label{eq:LQ.EL.syst.var}
\delta \ddot q &= A_1 \delta q + A_2 \delta \dot q + BR^{-1}B^\top \delta \kappa, \\
\delta \ddot \kappa &= \left(Q_2A_1-Q_1\right)\delta q + Q_2A_2\delta \dot q + \left( A_1^\top + Q_2BR^{-1}B^\top \right) \delta \kappa - A_2^\top \delta \dot \kappa.
\end{aligned}
\end{equation}
It is easy to see that \eqref{eq:LQ.Ham.syst.var} can be transformed into \eqref{eq:LQ.EL.syst.var} as follows. By taking the derivative with respect to time of equations on $\delta \dot q, \delta \dot \lambda_v$, it is possible to express \eqref{eq:LQ.Ham.syst.var} as a system of second order differential equations in $\delta q, \delta\lambda_v$. The transformation $\delta \lambda_v = \delta \kappa$ leads to \eqref{eq:LQ.EL.syst.var}. This shows that \eqref{eq:LQ.Ham.syst.var} and \eqref{eq:LQ.EL.syst.var} are equivalent and their solutions can be transformed one into another by $\delta \lambda_v = \delta \kappa$ and $\delta \lambda_q = Q_2 \delta \dot q - A_2^\top \delta \kappa - \delta \dot \kappa$. 
The characterisation of conjugate times of the LQ problem \eqref{LQ} was done in \cite{Agrachev2015} showing the dependence of the conjugate times on the eigendecomposition of the Hamiltonian matrix governing the linear system \eqref{eq:LQ.Ham.syst}. %{To check if smth can be said in case of the 2nd order control system or not}

\section{Conclusions and outlook}
{In this work, we extended the new Lagrangian formalism to solve OCPs \cite{sato2025} to include necessary and sufficient second order optimality conditions. We find that it is possible to define {corresponding} %the (strong)
Legendre conditions in the new framework and this formulation is also equivalent to the standard formulation due to the equivalence of the corresponding spaces given by an extended Tulzcyjew isomorphism. The new formulation allows {us to determine} %for determining the
conjugate points via the solution of {the} linearised Euler-Lagrange equation. The procedure of explicitly calculating the first conjugate time is given. Finally, we consider an example for which we derive the equations governing the conjugate time when the strong Legendre condition is satisfied.} 
%{some outlook}

The central advantage of the new Lagrangian formalism developed in \cite{sato2025} {lies} in the possibility to derive the numerical solutions of OCPs with second order dynamical constraints based on the discretisation of the new Lagrangian formulation as it has been done in \cite{kono2025}. The new framework of the second order conditions can be also transferred to the discrete setting, which will allow us to check the second order conditions numerically. Analysis of the second order conditions in the discrete setting and their relation to the conditions in the continuous setting will be the subject of future work. 
%%\begin{acknowledgement}
\section*{Acknowledgements}
The authors acknowledge the support of Deutsche Forschungsgemeinschaft (DFG) with the projects: LE
1841/12-1, AOBJ: 692092 and OB 368/5-1, AOBJ: 692093.
%%\end{acknowledgement}

\vspace{\baselineskip}
%% The style of the following references should be used in all documents.

%\bibliographystyle{ieee}
%\bibliography{references}
\printbibliography
\end{document}